\documentclass[12pt,twoside]{article}
\usepackage[english]{babel}
\usepackage{amsmath}
\usepackage{amssymb,amsfonts}
\usepackage{graphicx}                   

\begin{document}
\pagestyle{myheadings}
\markboth{\centerline{Jen\H o Szirmai}}
{The optimal hyperball packings ...}
\title
{The optimal hyperball packings related to the smallest compact arithmetic 5-orbifolds
}

\author{\normalsize Jen\H o Szirmai \\
\normalsize Budapest University of Technology and \\
\normalsize Economics Institute of Mathematics, \\
\normalsize Department of Geometry \\
\date{\normalsize{}}}


\maketitle


\begin{abstract}
The smallest three hyperbolic compact arithmetic 5-orbifolds
can be derived from two compact Coxeter polytops which are combinatorially simplicial prisms (or
complete orthoschemes of degree $d=1$) in the five dimensional hyperbolic space $\mathbf{H}^5$ (see \cite{BE} and \cite{EK}). The corresponding 
hyperbolic tilings are generated by reflections through their delimiting hyperplanes those involve the study of the relating
densest hyperball (hypersphere) packings with congruent hyperballs.

The analogous problem was discussed in \cite{Sz06-1} and \cite{Sz06-2} in the hyperbolic spaces $\mathbf{H}^n$ $(n=3,4)$.
In this paper we extend this procedure to determine the optimal hyperball packings to the above 5-dimensional
prism tilings. We compute their metric data and the
densities of their optimal hyperball packings, moreover, we formulate a conjecture for the candidate of the densest hyperball
packings in the 5-dimensional hyperbolic space $\mathbf{H}^5$.
\end{abstract}

\newtheorem{theorem}{Theorem}[section]
\newtheorem{conjecture}{Conjecture}[section]
\newtheorem{corollary}{Corollary}[section]
\newtheorem{lemma}{Lemma}[section]
\newtheorem{exmple}{Example}[section]
\newtheorem{defn}{Definition}[section]
\newtheorem{rmrk}{Remark}[section]
\newenvironment{definition}{\begin{defn}\normalfont}{\end{defn}}
\newenvironment{remark}{\begin{rmrk}\normalfont}{\end{rmrk}}
\newenvironment{example}{\begin{exmple}\normalfont}{\end{exmple}}
\newenvironment{acknowledgement}{Acknowledgement}



\section{Introduction}

Let in the 5-dimensional hyperbolic space be the group of isometries $Isom(\mathbf{H}^5)$ and
its orientationpreserving subgroup is denoted by $Isom^+(\mathbf{H}^5)$.
In \cite{BE} the lattice of smallest covolume
among cocompact arithmetic lattices of $Isom^+(\mathbf{H}^5)$ was determined.

In \cite{EK} the second and third values in the
volume spectrum of compact orientable arithmetic hyperbolic 5-orbifolds are determined and has proved the
following
\begin{theorem}[Emery-Kellerhals]
The lattices $\Gamma'_0$, $\Gamma'_1$, $\Gamma'_2$ (ordered by increasing covolume)
are the three cocompact arithmetic lattices in $Isom^+(\mathbf{H}^5)$ of minimal covolume.
They are unique, in the sense that any cocompact arithmetic lattice in
$Isom^+(\mathbf{H}^5)$ of covolume smaller than or equal to $\Gamma'_2$ is conjugate in $Isom(\mathbf{H}^5)$
to one of the $\Gamma'_i$ $(i=0,1,2)$.
\end{theorem}
The above lattices $\Gamma'_0$, $\Gamma'_1$, $\Gamma'_2$
can be derived by compact Coxeter polytopes $\mathcal{S}_i$ $(i=0,1,2)$ in $\mathbf{H}^5$ that are
characterized by the Coxeter symbols $[5,3,3,3,3]$, $[5,3,3,3,3^{1,1}]$ $[5,3,3,3,4]$ and Coxeter diagrams as follows:
\begin{figure}[ht]
\centering
\includegraphics[width=8cm]{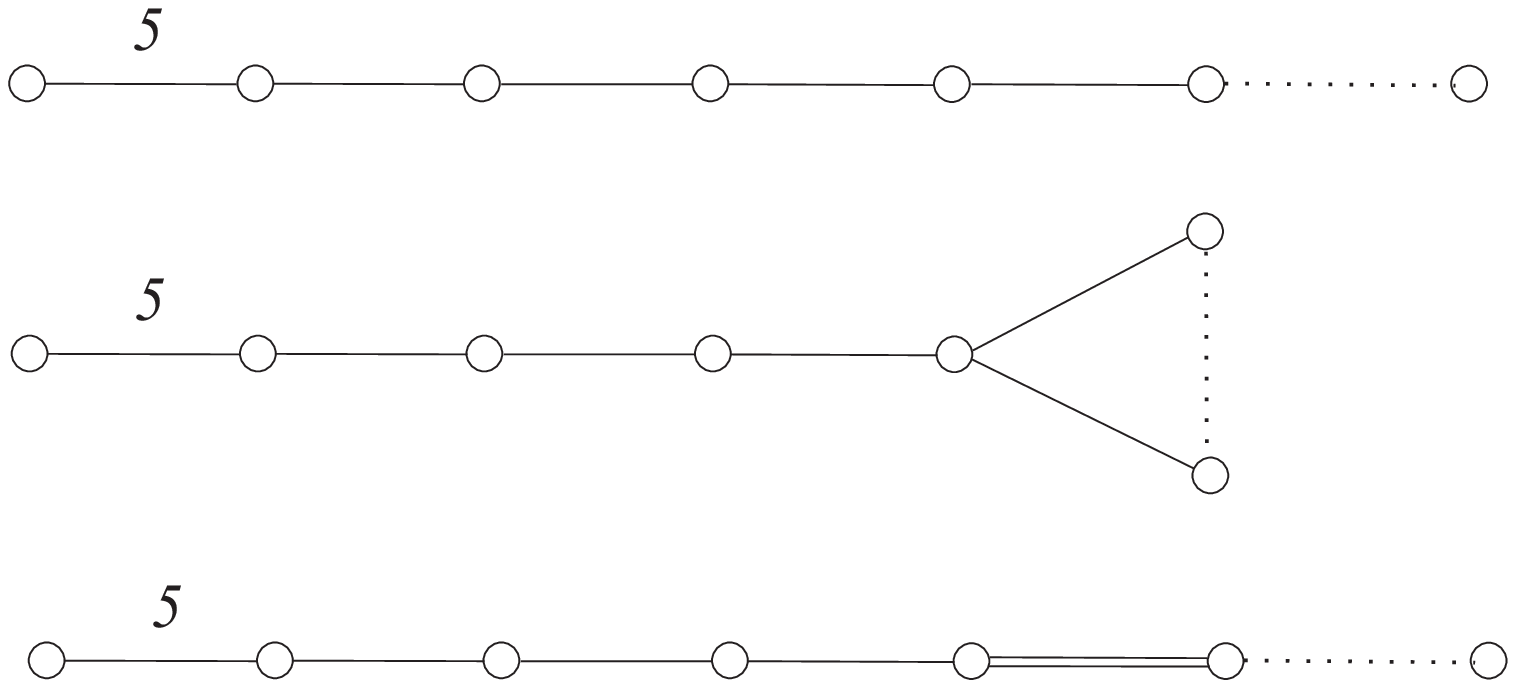}
\caption{}
\label{}
\end{figure}

The Coxeter group generated by the reflections
through the hyperplanes delimiting $\mathcal{S}_i$ $(i=0,1,2)$ denoted by $\Delta_i$. The next theorem has been proved in \cite{EK}:
\begin{theorem}[Emery-Kellerhals]
For $i = 0,1,2$ let $\Delta^+_i$ be the lattice $\Delta_i \cap Isom^+(\mathbf{H}^5)$, which is
of index two in $\Delta_i$. Then $\Delta^+_i$
is conjugate to $\Gamma'_i$ in $Isom(\mathbf{H}^5)$. In particular,
$\Delta_0$ realizes the smallest covolume among the cocompact arithmetic lattices
in $Isom(\mathbf{H}^5)$.
\end{theorem}
The polytops $\mathcal{S}_0$ and $\mathcal{S}_2$ are complete Coxeter orthoschemes of degree $d=1$ (see next Section) that were classified by {Im Hof}
in \cite{IH85} and \cite{IH90}. Moreover, the polytop $\mathcal{S}_1$
can be derived by the reflection of $\mathcal{S}_0$ in its facet $[5,3,3,3]$. 

In the hyperbolic space $\mathbf{H}^n$ $(n \ge 3)$ a regular prism is the convex hull of two congruent
$(n-1)$ dimensional regular polyhedra in ultraparallel hyperplanes, (i.e. $(n-1)$-planes), related by ,,translation" along the line
joining their centres that is the common perpendicular of the two hyperplanes.
Each vertex of such a tiling is either proper point or every vertex lies on the absolute
quadric of $\mathbf{H}^n$, in this case the prism tiling is called fully asymptotic.
Thus the prism is a polyhedron having at each vertex one $(n-1)$-dimensional regular polytop
and some $(n-1)$-dimensional prisms, meeting at this vertex. \newline 
\indent From the definitions of the regular prism tilings and the complete orthoschemes of degree $d=1$ (see next Section) follows that
a prism tiling exists in the $n$-dimensional hyperbolic space $\mathbf{H}^n, \ \ n \geq 3$ if and only if
exists an appropriate complete Coxeter orthoscheme of degree $d=1$.
The formulas for the hyperbolic covolumes of the considered 5-dimensional Coxeter tilings are determined in \cite{EK} (see also the formula (4.3))
therefore, it is possible to compute the covolumes of the regular prisms and the densities of the corresponding hyperball packings. \newline
\indent In \cite{Sz06-1} and \cite{Sz06-2} we have studied the regular prism tilings and their optimal hyperball packings in $\mathbf{H}^n$ $(n=3,4)$. \newline
\indent In this paper we extend the in former papers developed method to 5-dimensional hyperbolic space and
construct to each above described interesting Coxeter tiling a regular prism tiling in the 5-dimensional hyperbolic space,
study the corresponding optimal hyperball packings, moreover, we determine their metric data and their densities (see Table 1-2).
In the hyperbolic plane $\mathbf{H}^2$ the universal upper bound of the hypercycle packing density is $\frac{3}{\pi}$ 
determined by I.~Vermes in \cite{V79} 
and recently, (to the author's best knowledge) the candidates for the densest hyperball (hypersphere) packings in the $3$ and $4$-dimensional hyperbolic space 
$\mathbf{H}^n$ $n\ge3,4$ are derived by the regular prism tilings which are studied in papers \cite{Sz06-1}, \cite{Sz06-2} and in the present paper.
\vspace{3mm}

\centerline{\vbox{
\halign{\strut\vrule~\hfil $#$ \hfil~\vrule
&\quad \hfil $#$ \hfil~\vrule
&\quad \hfil $#$ \hfil~\vrule
\cr
\noalign{\hrule}
\multispan3{\strut\vrule\hfill\bf Table 1 \hfill\vrule}%
\cr
\noalign{\hrule}
\noalign{\vskip2pt}
\noalign{\hrule}
n& \mathrm{Coxeter ~ symbol} & \mathrm{The~known ~ maximal~ density} \cr
\noalign{\hrule}
3& [7,3,3] & \approx 0.82251367 \cr
\noalign{\hrule}
4& [3,5,3,3] & \approx 0.57680322 \cr
\noalign{\hrule}
5& [5,3,3,3,3]  & \approx 0.50514481 \cr
\noalign{\hrule}}}}
\smallbreak
\section{The projective model and \\ the complete orthoschemes }
We use for $\mathbf{H}^n$ the projective model in the Lorentz space $\mathbf{E}^{1,n}$ of signature $(1,n)$,
i.e.~$\mathbf{E}^{1,n}$ denotes the real vector space $\mathbf{V}^{n+1}$ equipped with the bilinear
form of signature $(1,n)$
\begin{equation}
\langle ~ \mathbf{x},~\mathbf{y} \rangle = -x^0y^0+x^1y^1+ \dots + x^n y^n \tag{2.1}
\end{equation}
where the non-zero vectors
$$
\mathbf{x}=(x^0,x^1,\dots,x^n)\in\mathbf{V}^{n+1} \ \  \text{and} \ \ \mathbf{y}=(y^0,y^1,\dots,y^n)\in\mathbf{V}^{n+1},
$$
are determined up to real factors, for representing points of $\mathcal{P}^n(\mathbf{R})$. Then $\mathbf{H}^n$ can be interpreted
as the interior of the quadric
\begin{equation}
Q=\{[\mathbf{x}]\in\mathcal{P}^n | \langle ~ \mathbf{x},~\mathbf{x} \rangle =0 \}=:\partial \mathbf{H}^n \tag{2.2}
\end{equation}
in the real projective space $\mathcal{P}^n(\mathbf{V}^{n+1},
\mbox{\boldmath$V$}\!_{n+1})$.

The points of the boundary $\partial \mathbf{H}^n $ in $\mathcal{P}^n$
are called points at infinity of $\mathbf{H}^n $, the points lying outside $\partial \mathbf{H}^n $
are said
to be outer points of $\mathbf{H}^n $ relative to $Q$. Let $P([\mathbf{x}]) \in \mathcal{P}^n$, a point
$[\mathbf{y}] \in \mathcal{P}^n$ is said to be conjugate to $[\mathbf{x}]$ relative to $Q$ if
$\langle ~ \mathbf{x},~\mathbf{y} \rangle =0$ holds. The set of all points which are conjugate to $P([\mathbf{x}])$
form a projective (polar) hyperplane
\begin{equation}
pol(P):=\{[\mathbf{y}]\in\mathcal{P}^n | \langle ~ \mathbf{x},~\mathbf{y} \rangle =0 \}. \tag{2.3}
\end{equation}
Thus the quadric $Q$ (by the symmetric bilinear form or scalar product in (2.1)) induces a bijection
(linear polarity $\mathbf{V}^{n+1} \rightarrow
\mbox{\boldmath$V$}\!_{n+1})$)
from the points of $\mathcal{P}^n$
onto its hyperplanes.

The point $X [\bold{x}]$ and the hyperplane $\alpha [\mbox{\boldmath$a$}]$
are called incident if $\bold{x}\mbox{\boldmath$a$}=0$ i.e. the value of the linear form
$\mbox{\boldmath$a$}$ on the vector $\bold{x}$
is equal to zero ($\bold{x} \in \bold{V^{n+1}} \setminus \{\mathbf{0}\}, \ \mbox{\boldmath$a$} \in \mbox{\boldmath$V$}_{n+1}
\setminus \{\mbox{\boldmath$0$}\}$).
The straight lines of $\mathcal{P}^n$ are characterized by 2-subspaces of $\bold{V^{n+1}} \ \text{or by $n-1$-spaces of} \
\mbox{\boldmath$V$}\!_{n+1}$, i.e. by 2 points or dually by $n-1$ hyperplane, respectively \cite{M97}.

Let $P \subset \mathbf{H}^n $denote a convex polytope bounded by finitely many hyperplanes $H^i$, which are
characterized by unit normal vectors $\mbox{\boldmath$b$}^i \in \mbox{\boldmath$V$}\!_{n+1}$ directed inwards with respect to $P$:
\begin{equation}
H^i:=\{\mathbf{x} \in \mathbf{H}^n | \langle ~ \mathbf{x},~\mbox{\boldmath$b$}^i \rangle =0 \} \ \ \text{with} \ \
\langle \mbox{\boldmath$b$}^i,\mbox{\boldmath$b$}^i \rangle = 1. \tag{2.4}
\end{equation}
We always assume that $P$ is acute-angled and of finite volume.

The Gram matrix $G(P):=( \langle \mbox{\boldmath$b$}^i, \mbox{\boldmath$b$}^j \rangle ) ~ {i,j \in \{ 0,1,2 \dots n \} }$
of the normal vectors $\mbox{\boldmath$b$}^i$ associated to $P$ is an indecomposable symmetric matrix of signature $(1,n)$
with entries $\langle \mbox{\boldmath$b$}^i,\mbox{\boldmath$b$}^i \rangle = 1$ and
$\langle \mbox{\boldmath$b$}^i,\mbox{\boldmath$b$}^j \rangle \leq 0$ for $i \ne j$, having the following
geometrical meaning
$$
\langle \mbox{\boldmath$b$}^i,\mbox{\boldmath$b$}^j \rangle =
\left\{
\begin{aligned}
&0 & &\text{if}~H^i \perp H^j,\\
&-\cos{\alpha^{ij}} & &\text{if}~H^i,H^j ~ \text{intersect \ on $P$ \ at \ angle} \ \alpha^{ij}, \\
&-1 & &\text{if}~\ H^i,H^j ~ \text{are parallel}, \\
&-\cosh{l^{ij}} & &\text{if}~H^i,H^j ~ \text{admit a common perpendicular of length} \ l^{ij}.
\end{aligned}
\right.
$$
A scheme $\varSigma$ is a weighted graph whose nodes $n_i,~n_j$ are joined by an edge with positive weight $\sigma^{ij}$
or are not joined at all; the last fact will be indicated by $\sigma^{ij}=0$. The number $\left| \varSigma \right|$
of nodes is called the order of $\varSigma$. To every scheme of order $m$ corresponds a symmetry matrix
$M(\varSigma)=(b^{ij})$ of order $m$ with $b^{ii}=1$ in the diagonal and non-positive entries $b^{ij}=-\sigma^{ij} \leq
0, \ i \ne j$, of it. The scheme $\varSigma(P)$ of an acute angled polytope $P$ is the scheme whose matrix
$M(\varSigma)$ coincides with the Gram matrix $G(P)$.

\begin{definition}
An orthoscheme $\mathcal{S}$ in $\mathbf{H}^n$ $(2\le n \in \mathbf{N})$ is a simplex bounded by $n+1$ hyperplanes $H^0,\dots,H^n$
such that
(see \cite{K91, B--H})
$$
H^i \bot H^j, \  \text{for} \ j\ne i-1,i,i+1.
$$
\end{definition}
\begin{rmrk}
This definition will be equivalent with the Definition 2.2:
\end{rmrk}
\begin{definition}
A simplex $\mathcal{S}$ in $\mathbf{H}^n$ is a orthoscheme iff the $n+1$ vertices of $\mathcal{S}$ can be
labelled by $A_0,A_1,\dots,A_n$ in such a way that
$$
\text{span}(A_0,\dots,A_i) \perp \text{span}(A_i,\dots,A_n) \ \ \text{for} \ \ 0<i<n-1.
$$
\end{definition}

Here we indicated the subspaces spanned by the corresponding vertices.

{\it The orthoschemes of degree} $d$ in $\mathbf{H}^n$ are bounded by $n+d+1$ hyperplanes
$H^0,H^1,\dots,H^{n+d}$ such that $H^i \perp H^j$ for $j \ne i-1,~i,~i+1$, where, for $d=2$,
indices are taken modulo $n+3$. For a usual orthoscheme we denote the $(n+1)$-hyperface opposite to the vertex $A_i$
by $H^i$ $(0 \le i \le n)$. An orthoscheme $\mathcal{S}$ has $n$ dihedral angles which
are not right angles. Let $\alpha^{ij}$ denote the dihedral angle of $\mathcal{S}$
between the faces $H^i$ and $H^j$. Then we have
\begin{equation}
\alpha^{ij}=\frac{\pi}{2}, \ \ \text{if} \ \ 0 \le i < j -1 \le n. \notag
\end{equation}
The $n$ remaining dihedral angles $\alpha^{i,i+1}, \ (0 \le i \le n-1)$ are called the
essential angles of $\mathcal{S}$.
Geometrically, complete orthoschemes of degree $d$ can be described as follows:
\begin{enumerate}
\item
For $d=0$, they coincide with the class of classical orthoschemes introduced by
{{Schl\"afli}} (see Definitions 2.1 and 2.3).
The initial and final vertices, $A_0$ and $A_n$ of the orthogonal edge-path
$A_iA_{i+1},~ i=0,\dots,n-1$, are called principal vertices of the orthoscheme (see Definition 2.3).
\item
A complete orthoscheme of degree $d=1$ can be interpreted as an 
orthoscheme with one outer principal vertex, say $A_n$, which is truncated by
its polar plane $pol(A_n)$ (see Fig.~2-3). In this case the orthoscheme is called simply truncated with
outer vertex $A_n$.
\item
A complete orthoscheme of degree $d=2$ can be interpreted as an 
orthoscheme with two outer principal vertices, $A_0,~A_n$, which is truncated by
its polar hyperplanes $pol(A_0)$ and $pol(A_n)$. In this case the orthoscheme is called doubly
truncated. (In this case we distinguish two different types of orthoschemes but I
will not enter into the details (see \cite{K89}, \cite{K91}).)
\end{enumerate}
For the schemes of complete Coxeter orthoschemes $\mathcal{S} \subset \mathbf{H}^n$ we adopt the usual conventions and use
them sometimes even in the Coxeter case: If two nodes are related by the weight $\cos{\frac{\pi}{p}}$
then they are joined by a ($p-2$)-fold line for $p=3,~4$ and by a single line marked $p$ for $p \geq 5$.
In the hyperbolic case if two bounding hyperplanes of $S$ are parallel, then the corresponding nodes
are joined by a line marked $\infty$. If they are divergent then their nodes are joined by a dotted line.

The principal minor matrix $(c^{ij})$ of $G(\mathcal{S})$ is the so called Coxeter-Schl\"afli matrix of the orthoschem $S$ with
parameters $p,~q,~r,~s,~t$:
\[
(c^{ij}):=\begin{pmatrix}
1& -\cos{\frac{\pi}{p}}& 0 & 0 & 0 & 0\\
-\cos{\frac{\pi}{p}} & 1 & -\cos{\frac{\pi}{q}}& 0 & 0 &0 \\
0 & -\cos{\frac{\pi}{q}} & 1 & -\cos{\frac{\pi}{r}} & 0 &0 \\
0 & 0 & -\cos{\frac{\pi}{r}} & 1 & -\cos{\frac{\pi}{s}} &0 \\
0 & 0 & 0 & -\cos{\frac{\pi}{s}} & 1&-\cos{\frac{\pi}{t}} \\
0& 0 & 0 & 0 & -\cos{\frac{\pi}{t}} & 1 \\
\end{pmatrix}. \tag{2.5}
\]
\begin{figure}[ht]
\centering
\includegraphics[width=8cm]{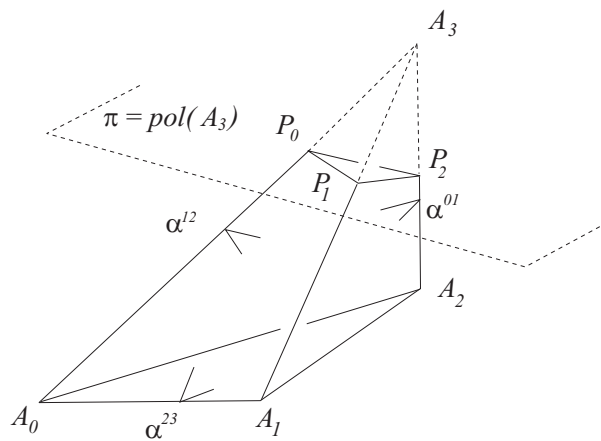}
\caption{Complete orthoscheme of degree $d=1$}
\label{}
\end{figure}
\section{Regular prism tilings and their optimal hyperball packings in $\mathbf{H}^5$}
\subsection{The structure of the 5-dimensional regular prism tilings}
In hyperbolic space $\mathbf{H}^n$ $(n \ge 3)$ a regular prism is the convex hull of two congruent
$(n-1)$ dimensional regular polyhedra in ultraparallel hyperplanes, (i.e. $(n-1)$-planes), related by ,,translation" along the line
joining their centres that is the common perpendicular of the two hyperplanes. The two regular 4-faces of a regular prism are called cover-polytops, 
and its other 4-dimensional facets are called side-prisms.

In this section we consider the 5-dimensional regular prism tilings $\mathcal{T}_{i}$ $(i=0,2)$
derived by Coxeter tilings $[5,3,3,3,3]$ and $[5,3,3,3,4]$.
(From the Coxeter tiling $[5,3,3,3,3^{1,1}]$ generated regular prism tiling is congruent to that derived by $[5,3,3,3,3]$.)

Fig.~3 shows a part of our 5-prism where $A_4$ is the centre of a cover-polyhedron,
$A_3$ is the centre of a 3-face of the cover polyhedron, $A_2$ is the midpoint of its $2$-face,
$A_1$ is a midpoint of an edge of this face, and $A_0$ is one vertex (end) of that edge.
\begin{figure}[ht]
\centering
\includegraphics[width=9cm]{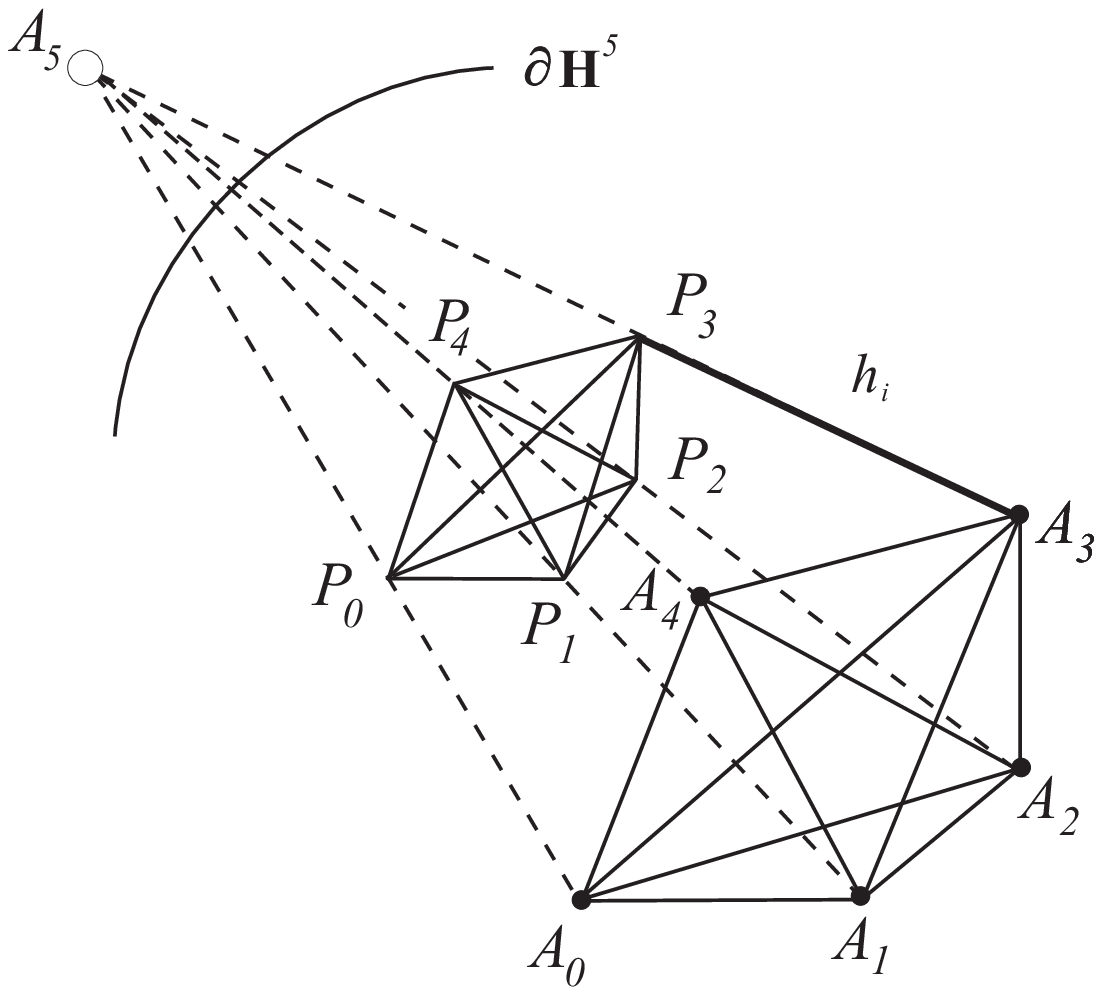}
\caption{}
\label{}
\end{figure}
Let $B_0,~B_1,~B_2,~B_3,B_4$ be the corresponding points of the other cover-polytop of the regular 5-prism.
The midpoints of the edges which do not lie in the cover-polytops
form a hyperplane denoted by $\pi$.
The foot points
$P_i ~ (i \in \{ 0,1,2,3,4 \})$
of the perpendiculars dropped from the points
$A_i$ on the plane $\pi$ form in both cases the same {\it{characteristic (or fundamental) simplex}} $\mathcal{S_{\pi}}$
with Coxeter-Schl\"afli symbol $[5,3,3,3]$ (see Fig.~3).

\begin{rmrk}
In $\mathbf{H}^3$ (see \cite{Sz06-1}) the corresponding prisms are called regular $p$-gonal prisms $(p \ge 3)$ in
which the regular polyhedra (the cover-faces)
are regular $p$-gons, and the side-faces are rectangles.
Fig.~4 shows a part of such a prism where $A_2$ is the centre of a regular $p$-gonal
face, $A_1$ is a midpoint of a side of this face, and $A_0$ is one vertex (end) of that side.
Let $B_0,~B_1,~B_2$ be the corresponding points of the other $p$-gonal face of the prism.
\end{rmrk}

Analogously to the 3-dimensional case, it can be seen that $\mathcal{S}_{i}=A_0A_1A_2A_3$ $A_4A_5$ $P_0P_1P_2P_3P_4P_5$ $(i=0,2)$ is an complete
orthoscheme with degree $d=1$ where $A_5$ is an outer vertex of
$\mathbf{H}^5$ and the points $P_0,P_1,P_2,P_3,P_4$ lie in its polar hyperplane $\pi$ (see Fig.~3).
The corresponding regular prism $\mathcal{P}_{i}$ can be derived by reflections in facets of $\mathcal{S}_{i}$ containing the
point $A_3$.

{\it We consider the images of $\mathcal{P}_{i}$ under reflections in its side facets (side prisms).
The union of these 5-dimensional regular prisms (having the common $\pi$ hyperplane) forms an infinite polyhedron denoted by $\mathcal{F}_{i}$.}

From the definitions of the prism tilings and the complete orthoschemes of degree $d=1$ follows that a
regular prism tiling exists in the $n$-dimensional hyperbolic space $\mathbf{H}^n, \ \ n \geq 3$ if and only if
exists a complete Coxeter orthoscheme of degree $d=1$ with two divergent faces.

The complete Coxeter orthoschemes were classified by {{Im Hof}}
\cite{IH85} and \cite{IH90} by generalizing the method of {{Coxeter}} and {{B\"ohm}} appropriately.
He showed that they exist only for dimensions $\leq 9$.

On the other hand, if a 5-dimensional regular prism tiling $[p,q,r,s,t]$ exists, then it has to satisfy the following two requirements:
\begin{enumerate}
\item The orthogonal projection of the cover-polytop
on the hyperbolic hyperplane $\pi$ is a regular Coxeter honeycomb with proper vertices and centres.
Using the classical notation of the tesselations, these honeycombs are given by their
Coxeter-Schl\"afli symbols $[ p,q,r,s ]$.
\item The vertex figures about a vertex of
such a prism tiling has to form a 5-dimensional regular polyhedron.
\end{enumerate}
\begin{figure}[ht]
\centering
\includegraphics[width=10cm]{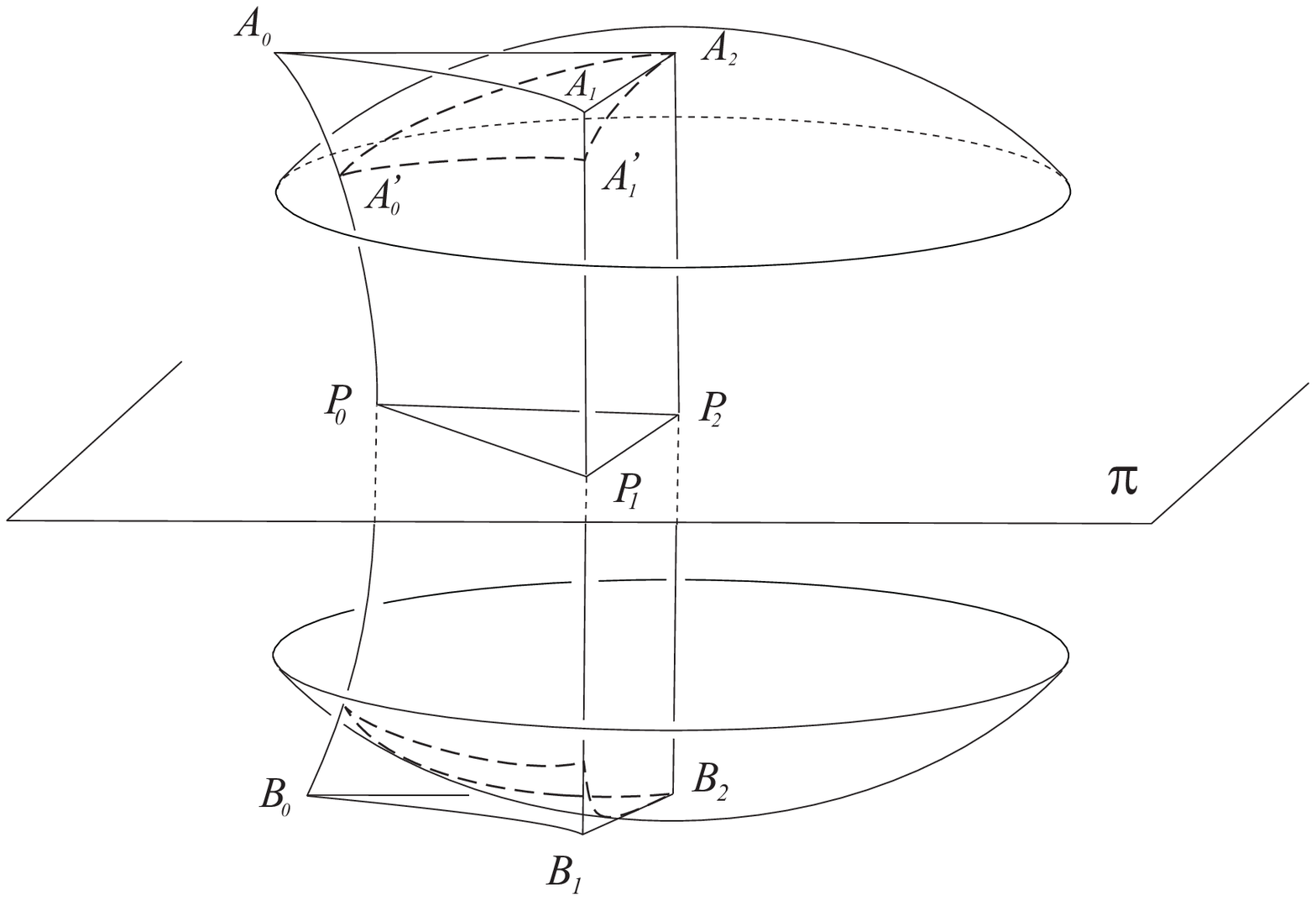}
\caption{}
\label{}
\end{figure}
\subsection{The volumes of the orthoschemes}
A plane orthoscheme is a right-angled triangle, whose area formula can be expressed by the
well known defect formula. For three-dimensional spherical orthoschemes, {{L.~Schl\"afli}} about 1850
was able to find the volume differentials depending on differential of the not fixed
3 dihedral angles. Already in 1836, {{N.~I.~Lobachevsky}} found
a volume formula for three-dimensional hyperbolic orthoschemes $\mathcal{S}$ \cite{B--H}.

The integration method for orthoschemes of dimension three was generalized by
{{B\"ohm}} in 1962 \cite{B--H} to spaces of constant nonvanishing curvature of arbitrary
dimension.

{{R.~Kellerhals}} in \cite{K89} derived a volume formula in the 3 dimensional hyperbolic space for the
{\it{complete orthoschemes of degree}} $d,~(d=0,~1,~2)$ and she explicitly determined in \cite{K91} the volumes of
all complete hyperbolic orthoschemes in even dimension ($n \geq 4$),
moreover, she in \cite{K92} developed a procedure to determine the volumes of 5-dimensionalal hyperbolic orthschemes.

The volumes of the complete orthoschemes $\mathcal{S}_i$ $(i=0,2)$ of $d=1$ can be computed by volume differential formula of L.~ Schl\"afli
with the the following formula (see \cite{EK}):
\begin{equation}
Vol_5(\mathcal{S}_i)=\frac{1}{4}\int_{\alpha_i}^{\frac{2\pi}{5}}{Vol_3([5,3,\beta(t)])dt+\frac{\zeta(3)}{3200}} \tag{3.1}
\end{equation}
with a compact tetrahedron $[5,3,\beta(t)]$ whose angle parameter $0<\beta(t)<\frac{\pi}{2}$ is given by
$$
\beta({t})=\arctan{\sqrt{2-\cot^2{t}}}.
$$
Then, the volume of the 3-dimensional orthoscheme face $[5,3,\beta(t)]$ as given by Lobachevsky's formula:
\begin{equation}
\begin{gathered}
Vol_3([5,3,\beta(t)])=\frac{1}{4} \{\mathcal{L}_2\big(\frac{\pi}{5}+\theta(t)\big)-\mathcal{L}_2\big(\frac{\pi}{5}-\theta(t)\big)-\mathcal{L}_2\big(\frac{\pi}{6}+\theta(t)\big)+\\
\mathcal{L}_2\big(\frac{\pi}{6}-\theta(t)\big)+\mathcal{L}_2\big(\beta(t)+\theta(t)\big)-\mathcal{L}_2\big(\beta(t)-\theta(t)\big)+2\mathcal{L}_2\big(\frac{\pi}{2}-\theta(t)\big) \end{gathered} \tag{3.2}
\end{equation}
where
$$
\mathcal{L}(\omega)=-\int_0^{\omega}{\log|2\sin{t}|dt}, \ \omega \in \mathbf{R},
$$
is the Lobachevsky's function and
$$
\theta(t)=\arctan\frac{\sqrt{1-4\sin^2\frac{\pi}{5} \sin^2{\beta(t)}}}{2\cos\frac{\pi}{5} \cos\beta(t)}.
$$
\subsection{The optimal hyperball packing}
The equidistance surface (or hypersphere) is a quadratic surface at a constant distance
from a plane in both halfspaces. The infinite body of the hypersphere is called hyperball.

The 5-dimensional hypersphere with distance $h$ to a hyperplane $\pi$ is denoted by $\mathcal{H}^h$. The volume of a bounded hyperball piece 
$\mathcal{H}^h(\mathcal{A})$ delimited by the 4-polytop $\mathcal{A} \subset \pi$, $\mathcal{H}^h$ and some to $\pi$ orthogonal 
hyperplanes derived by the facets of $\mathcal{A}$
can be determined by the formula (3.3) that follows by the generalization of the classical method of {{J.~Bolyai}}:
\begin{equation}
Vol({\mathcal{H}}^h(\mathcal{A}))=\frac{1}{16} Vol(\mathcal{A})k  \Big( \frac{1}{2} \sinh \frac{4h}{k}+
4 \sinh \frac{2h}{k}\Big) +\frac{3h Vol(\mathcal{A})}{8}, \tag{3.3}
\end{equation}
where the volume of the hyperbolic 4-polytop $\mathcal{A}$ lying in the plane $\pi$ is $Vol(\mathcal{A})$.
The constant $k =\sqrt{\frac{-1}{K}}$ is the natural length unit in
$\mathbf{H}^n$. $K$ will be the constant negative sectional curvature.

{\it{We are looking for the optimal hyperball $\mathcal{H}^{h_i}_{opt}$ inscribed in $\mathcal{T}_i$ with maximal height.}}

The optimal hypersphere $\mathcal{H}_{opt}^{h_i}$ touches the cover-faces of the regular 5-prisms containing by $\mathcal{F}_i$.
Therefore, the optimal distance from the 4-midplane $\pi$
will be $h_{i}=P_4A_4$ (Fig.~3).

We consider one from the former defined 5-dimensional regular prism tilings $\mathcal{T}_i$ $(i=0,2)$ and
the infinite polyhedron $\mathcal{F}_i$ derived from that (the union of 5-dimensional regular prisms having the common hyperplane $\pi$).
$\mathcal{F}_i$ and its images under reflections in its "cover facets" fill the hyperbolic
space $\mathbf{H}^5$ without overlap thus we obtain by the above images of $\mathcal{H}_{opt}^{h_i}$
a locally optimal hyperball packing to the tiling $\mathcal{T}_{i}$ $(i=0,2)$.

The points $P_4({\mathbf{p}}_4)$ and $A_4({\mathbf{a}}_4)$ are proper points of the hyperbolic 5-space and
$P_4$ lies on the polar hyperplane $pol(A_5)$ of the outer point $A_5$ thus
\begin{equation}
\begin{gathered}
\mathbf{p}_4=c \cdot \mathbf{a}_5+\mathbf{a}_4 \in \mbox{\boldmath$a$}^5 \Leftrightarrow
c \cdot \mathbf{a}_5 \mbox{\boldmath$a$}^5+\mathbf{a}_4 \mbox{\boldmath$a$}^5=0 \Leftrightarrow
c=-\frac{\mathbf{a}_4 \mbox{\boldmath$a$}^5}{\mathbf{a}_5 \mbox{\boldmath$a$}^5} \Leftrightarrow \\
\mathbf{p}_4 \sim -\frac{\mathbf{a}_4 \mbox{\boldmath$a$}^5}{\mathbf{a}_5 \mbox{\boldmath$a$}^5}
\mathbf{a}_5+\mathbf{a}_4 \sim \mathbf{a}_4 (\mathbf{a}_5 \mbox{\boldmath$a$}^5) - \mathbf{a}_5 (\mathbf{a}_4 \mbox{\boldmath$a$}^5)=
\mathbf{a}_4 h_{55}-\mathbf{a}_5 h_{45},
\end{gathered} \tag{3.4}
\end{equation}
where $h_{ij}$ is the inverse of the Coxeter-Schl\"afli matrix $c^{ij}$ (see (2.5)) of the orthoscheme $\mathcal{S}_i$.
The hyperbolic distance $h_i$ can be calculated by the following formula \cite{M89}:
\[
\begin{gathered}
\cosh{P_4A_4}=\cosh{h_i}=\frac{- \langle {\mathbf{p}}_4, {\mathbf{a}}_4 \rangle }
{\sqrt{\langle {\mathbf{p}}_4, {\mathbf{p}}_4 \rangle \langle {\mathbf{a}}_4, {\mathbf{a}}_4 \rangle}}= \\ =\frac{h_{45}^2-h_{44}h_{55}}
{\sqrt{h_{44}\langle \mathbf{p}_4, \mathbf{p}_4 \rangle}} =
\sqrt{\frac{h_{44}~h_{55}-h_{45}^2}
{h_{44}~h_{55}}}.
\end{gathered} \tag{3.5}
\]
The volume of the polyhedron $\mathcal{S}_i$ is denoted by $Vol_5(\mathcal{S}_i)$ (see Section 2).

For the density of the packing it is sufficient to relate the volume of the optimal hyperball piece to that of
its containing polyhedron $\mathcal{S}_i$
(see Fig.~3) because the tiling can be constructed of such polyhedron. This polytope
and its images in $\mathcal{F}_i$ divide the $\mathcal{H}^{h_i}$ into congruent
pieces whose volume is denoted by $Vol({\mathcal{H}_{opt}^{h_i}})$. We illustrate in the 3-dimensional case such a
hyperball piece $A_2 A_0' A_1' P_0 P_1 P_2$ in Fig.~4.

The density of the optimal hyperball
packing to the prism tiling $\mathcal{T}_{i}$ $(i=0,2)$ is defined by the following formula:

\begin{definition}
\begin{equation}
\delta^{opt}(\mathcal{T}_{i}):=\frac{Vol(\mathcal{H}_{opt}^{h_i})}{Vol_5({\mathcal{S}_i})}. \tag{3.6}
\end{equation}
\end{definition}
$\delta^{opt}(\mathcal{T}_{i})$ can be determined by the formulas (3.1), (3.2), (3.3) and (3.5) using, that the volume $Vol(\mathcal{A})=Vol(\mathcal{S}_{\pi})$ (see (3.3)) of the 4-dimensional polytope $\mathcal{S}_{\pi}$ with Coxeter symbol $[5,3,3,3]$ in both cases is
$$
Vol(\mathcal{S}_{\pi})=\frac{\pi^2}{10800} \approx0.00091385 .
$$

Finally we get the following results:
\medbreak
\centerline{\vbox{
\halign{\strut\vrule~\hfil $#$ \hfil~\vrule
&\quad \hfil $#$ \hfil~\vrule
&\quad \hfil $#$ \hfil~\vrule
\cr
\noalign{\hrule}
\multispan3{\strut\vrule\hfill\bf Table 1 \hfill\vrule}%
\cr
\noalign{\hrule}
\noalign{\vskip2pt}
\noalign{\hrule}
&\varSigma_{53333} & \varSigma_{53334} \cr
\noalign{\hrule}
Vol(S_i)& \approx 0.00076730 & \approx 0.00198469 \cr
\noalign{\hrule}
h_i& \approx 0.38359861  & \approx 0.53063753 \cr
\noalign{\hrule}
Vol(\mathcal{H}^{h_i}_{opt})& \approx 0.00038760 & \approx 0.00059001 \cr
\noalign{\hrule}
\delta^{opt}& \approx 0.50514481 & \approx 0.29727979 \cr
\noalign{\hrule}}}}
\smallbreak
\begin{rmrk}
The optimal density of the Coxeter tiling $\mathcal{T}_1$ with Coxeter symbol $[5,3,3,3,3^{1,1}]$ is equal to $\delta^{opt}(\mathcal{T}_0) 
\approx 0.50514481$.
\end{rmrk}
The next conjecture for the optimal hyperball packings can be formulated.
\begin{conjecture}
{The above described optimal hyperball packings to Coxeter tilings $[5,3,3,3,3]$ and $[5,3,3,3,3^{1,1}]$ provide the densest hyperball packings 
in the 5-dimensional hyperbolic space $\mathbf{H}^5$.}  
\end{conjecture}
\begin{rmrk}
Regular hyperbolic honeycombs exist only up to 5 dimensions \cite{C56}. Therefore regular prism tilings
can exist up to 6 dimensions. From the definitions of the prism tilings and the complete orthoschemes of degree
$d=1$ it follows that prism tilings exist in the $n$-dimensional hyperbolic space $\mathbb{H}^n, \ \ n \geq 3$ if and only if
there exist complete Coxeter orthoschemes of degree $d=1$ with two divergent faces.
From the paper \cite{IH90} it follows that in the 5-dimensional hyperbolic space there is a further type and in the 6-dimensional hyperbolic space there is no such a Coxeter orthoschem.
\end{rmrk}


%
\normalsize Jen\H o Szirmai, ~ 
\normalsize Budapest University of Technology and \\
\normalsize Economics, Institute of Mathematics, \normalsize Department of Geometry \\
\normalsize H-1521 Budapest, Hungary \\
\normalsize Email: ~ szirmai@math.bme.hu,~ 
\normalsize www.math.bme.hu /$\sim$szirmai
\end{document}